\documentclass{amsart}%
\usepackage{graphicx}
\usepackage{amscd}
\usepackage{amsmath}
\usepackage{amsfonts}
\usepackage{amssymb}%
\newtheorem{theorem}{Theorem}
\theoremstyle{plain}

\numberwithin{equation}{section}

\thanks
{Y.Sezer}
\thanks
{Department Of Mathematics, Faculty of Arts and Science, Y\i ld\i z Technical
University,Davutpa\d{s}a, \.Istanbul, Turkey }
\thanks
{e-mail: {\tt ysezer@yildiz.edu.tr}}

\thanks
{\"{O}. Baksi}
\thanks
{Department Of Mathematics, Faculty of Arts and Science, Y\i ld\i z Technical
University,Davutpa\d{s}a, \.Istanbul, Turkey}
\thanks
{e-mail: {\tt baksi@yildiz.edu.tr}}

\thanks
{\emph{Mathematics Subject Classification:} 47A10, 34L20}

\begin{document}

\leftline {\qquad Yonca SEZER, \"{O}zlem BAKSI}
 
\vskip 0.8cm

\title {The Second Regularized Trace of Even Order Differential Operators with Operator}

\maketitle

\date{Received: date / Accepted: date}

\vskip 1.8cm

\noindent{\bf{Abstract.}} In this paper, we investigate the spectrum of the self adjoint operator $L$ defined by
$$L:=(-1)^{r}\frac{d^{2r}}{dx^{2r}}+A+Q(x),$$ where $A$ is a self adjoint operator and $Q(x)$ is a nuclear operator in a separable Hilbert
space. We also derive asymptotic formulas for the sum of eigenvalues of the operator $L$.

\vskip 1cm

\section {\bf{Introduction}}

\noindent The theory of regularized traces of differential operators began with the study of Gelfand and Levitan \cite{Ge}. They calculated
the trace formula for the sum of substraction of eigenvalues of two self adjoint operators. After this primary work, many mathematicians
concentrated on this theory in a large scale.

\noindent Dikiy \cite{Di}, Halberg and Kramer \cite{Ha}, Levitan \cite{Lv} and some others studied the regularized traces of scalar
differential operators. The list of works on the subject was given in the works Levitan and Sargsyan \cite{Le} and Fulton and Pruess
\cite{Fu}, but a few of these works were about the regularized trace of differential operators with operator coefficient. Chalilova \cite{Ch}
calculated regularized trace of Sturm Lioville operator with bounded operator coefficient. Ad\i g\"uzelov \cite{Ad} computed regularized
trace of the difference of two Sturm-Liouville operators with bounded operator coefficient given in the semi-axis. Maksudov, Bayramoglu and
Ad\i g\"uzelov \cite{Ma} found a formula for the regularized trace of  Sturm-Liouville operators with unbounded operator coefficient under
the Dirichlet boundary conditions.  Bayramoglu and Ad\i g\"uzelov \cite{Ba} obtained the regularized trace of second order singular
differential operator with bounded operator coefficient. Furthermore Ad\i g\"uzelov and Bak\d{s}i \cite{Oz},  Ad\i g\"uzelov and Sezer \cite{Yo}, \cite{Pr} investigated the regularized trace formulas of differential operator with operator coefficient.

\noindent Although most of the previous researches on the subject dealt with regularized trace of second order differential operators, we
focused on higher order differential operators. It is clear that our study advances the formulation of regularized trace that the prior manuscripts has proved. This paper aims to explore the second regularized trace of higher order differential operators with operator coefficient.

\noindent Let us begin by recalling some definitions and properties:

\noindent Let $\>H\>$ be an infinite dimensional separable Hilbert space.  We denote  the inner products in $\>H\>$ by  $\>(.,.)\>$ and the
norm in $\>H\>$ by $ \|.\|$. Let $H_{1}=L_{2}(0,\pi;\>H)$ denote the set of all functions $f$ from $[0,\pi]$ into $H$ which are strongly
measurable and satisfy the condition $\int_{0}^{\pi}\bigl\|f(x)\bigr\|^{2}dx<\infty$.The space $H_{1}$ is a linear space. If the inner
product of arbitrary two elements $f$ and $g$ of the space $H_{1}$ is defined as $(f,g)_{H_{1}}=\int_{0}^{\pi}\bigl(f(x),g(x)\bigr)dx$, then
$H_{1}$ becomes an infinite dimensional separable Hilbert space \cite{Kr}. The norm in the space $H_{1}$ is denoted by  $ \|.\|_{1}$.

\noindent $\sigma_{\infty}(H)$ denotes the set of all compact operators from $H$ to $H$. If $T \in \sigma_{\infty}(H)$, then $T^{*}T$ is a
nonnegative self adjoint operator and $(T^{*}T)^{\frac{1}{2}}\in\sigma_{\infty}(H) $  \cite{Co}. Let the nonzero eigenvalues of the operator
$(T^{*}T)^{\frac{1}{2}}$ be $ \{s_{j}\}_{j=1}^{k}$ $(0\leq k\leq\infty)$ such that $s_{1}\geq s_{2}\geq...\geq s_{k}$ according to its multiplicity number. Since $(T^{*}T)^{\frac{1}{2}}$ is non negative, $s_{k}$'s are positive numbers. The numbers
$s_{k}$ are called s-numbers of the operator $T$. If $k< \infty$, then $s_{j}=0$ $(j=k+1,k+2,...)$ will be accepted. s-numbers of the
operator $T$ are also denoted by $s_{j}(T)$ $(j=1,2,...)$. Here $s_{1}(T)=\|T\|$.

\noindent  If $T$ is a normal operator, then $s_{j}(T)=|\lambda_{j}(T)|$ $(j=1,2,...,k)$ \cite{Co}. Here,
$\{\lambda_{1}(T),\lambda_{2}(T),...,\lambda_{k}(T)\}$ is an ordering of all nonzero eigenvalues of the operator $T$ according to
$|\lambda_{1}(T)|\geq|\lambda_{2}(T)|\geq...\geq|\lambda_{k}(T)|$. $\sigma_{p}$ or $\sigma_{p}(H)$ denotes the set of all compact operators,
the s-numbers of which satisfy the condition $\sum\limits_{j=1}^{\infty}s_{j}^{p}(T)<\infty$ $(p\geq1)$. The set $\sigma_{p}$ $(p\geq1)$ is a
separable Banach space with respect to the norm $\|T\|_{\sigma_{p}(H)}=\Biggl[\sum\limits_{j=1}^{\infty}s_{j}^{p}(T)\Biggr]^{\frac{1}{p}}$,
for $T\in\sigma_{p}(H)$ \cite{Co}. $\sigma_{1}(H)$ is the set of all operators $T\in\sigma_{\infty}(H)$, the s-numbers of which satisfy the condition
$\sum\limits_{j=1}^{\infty}s_{j}(T)<\infty$ . If an operator belongs to $\sigma_{1}(H)$, then it is called a nuclear operator. If the
operators $T \in\sigma_{p}(H)$ and $B \in B(H)$, then $TB$,$BT \in \sigma_{p}(H)$ and $\|TB\|_{\sigma_{p}(H)}\leq\|B\|\|T\|_{\sigma_{p}(H)}$, $\|BT\|_{\sigma_{p}(H)}\leq\|B\|\|T\|_{\sigma_{p}(H)}.$

\noindent If $T$ is a nuclear operator and $\{e_{j}\}_{j=1}^{\infty}\subset H$ is any orthonormal basis, then the series
$\sum\limits_{j=1}^{\infty}\Bigl(Te_{j},e_{j}\Bigr)$ is convergent and the sum of the series  does not depend on the choice of the basis
$\{e_{j}\}_{j=1}^{\infty}$. The sum of this series is said to be matrix trace of the operator $T$ denoted by $trT$,
$$trT=\sum\limits_{j=1}^{\nu(T)}\lambda_{j}(T) \eqno(1.1)$$

\noindent \cite{Co}. Here, each eigenvalue counted according to its algebraic multiplicity number.  $\nu(T)$ denotes the sum of algebraic
multiplicity of non-zero eigenvalues of the operator $T$ \cite{Co}. The sum of the series $\sum\limits_{j=1}^{\nu(T)}\lambda_{j}(T)$ is
called spectral trace of the operator $T$.

\noindent Now, let us return to our problem.

\noindent Consider the differential expression
$$\ell_{0}(y)=(-1)^{r}y^{(2r)}(x)+Ay(x)\eqno(1.2)$$

\noindent  in the space $H_{1}=L_{2}(0,\pi;\>H)$. Here, the densely defined operator

\noindent $\>A:D(A) \rightarrow H\> $ satisfies the conditions $A=A^{*}\geq I $ ($I$ is unit operator in $H$) and

\noindent $A^{-1}\in\sigma_{\infty}(H)$. Let $\>\{\gamma_{n}\}_{1}^{\infty} \>$ be  an ordering of all eigenvalues of $A$ according to
 $\>\gamma_{1}\leq\gamma_{2}\leq\cdots\leq\gamma_{n}\leq\cdots \>$ and $\varphi_{n}$ the corresponding orthonormal eigenfunctions. Here, each eigenvalue counted according to its multiplicity
 number.

\noindent Let $\>D_{0}\>$ be a subset of the space $H_{1}\>$. A function $\>y(x)\in D_{0}\>\>$, if $y(x)$ satisfies the following conditions:

\noindent\textbf{(y1)} $\>y(x)\>$ has continuous derivative of the $(2r)$th-order with respect to the norm in the space  $\>H\>$ for every
$x\in\>[0,\pi]\>$,

\noindent\textbf{(y2)} $\>Ay(x)\>$ is continuous with respect to the norm of the space $\>H \>$ on $\>[0,\pi]\>$,

\noindent\textbf{(y3)} $\> y^{'}(0)=y^{'''}(0)= \cdots=y^{(2r-1)}(0)=y(\pi)=y^{''}(\pi)= \cdots =y^{(2r-2)}(\pi)=0\>$ $(r=1,2,\dots,m)$.

\noindent Here, $\>\overline{D_{0}}=H_{1}\>$. Define the linear operator $\>L^{'}_{0}:D_{0}\rightarrow\>H_{1}\>$ as $L^{'}_{0}y:=\ell_{0}(y)$.

\noindent The construction above gives that $L^{'}_{0}$  is symmetric. The eigenvalues of $L^{'}_{0} $ are

\noindent $\Bigl(k+\frac{1}{2}\Bigr)^{2r}+\gamma_{j}$ $(k=0,1,2,\ldots;j=1,2,\ldots)$ and $\sqrt{\frac{2}{\pi}}\varphi_{j}\cos\Bigl(k+\frac{1}{2}\Bigr)x $ the corresponding orthonormal eigenvectors.

\noindent We can see that the orthonormal eigenvector system of the symmetric operator $L^{'}_{0} $ is an orthonormal basis in the space
$H_{1}$. We denote the closure of $L^{'}_{0} $ by $\>L_{0}$ defined as $\>L_{0}:D(L_{0})\rightarrow H_{1}\>$. Since the orthonormal
eigenvector system of the operator $L^{'}_{0} $ is an orthonormal basis in the space $H_{1}$,  $\>L_{0}$  is a self adjoint operator.

\noindent Let $\>Q(x)\>$ defined on $ [0,\pi]$ be an operator function satisfying the following conditions:

\noindent\textbf{(Q1)} $\>Q(x)\>$ has weak derivative of $(2r+2)$th order  and

\noindent $Q^{(2i+1)}(0)+Q^{(2i+1)}(\pi)=0 \qquad (i=0,1,2,\ldots,r)$

\noindent\textbf{(Q2)} $\>Q^{(i)}(x):H\rightarrow H\qquad (i=0,1,2,\ldots,2r+2) \>$ are self-adjoint operators for every $\>x \in [0,\pi]\>$,
$AQ''(x),Q^{(2r+2)}(x)\in\sigma_{1}(H)$ and the functions $\>\|AQ''(x)\|_{\sigma_{1}(H)}$, $\>\|Q^{(2r+2)}(x)\|_{\sigma_{1}(H)}$ are bounded
and measurable in the interval $\>[0,\pi]\>$.

\noindent The operator $L: D(L_{0})\rightarrow H_{1}$ defined by
$$L=L_{0}+Q$$
is a self adjoint operator.  The operators $L_{0}$ and $L$
have purely discrete spectrum \cite{Oz}. We denote the resolvent sets of $L_{0}$, $L$ by $\rho (L_{0})$, $\rho (L)$ and the resolvent operators of $L_{0}$, $L$ by $R_{\lambda}^{0}=(L_{0}-\lambda I)^{-1}$ , $R_{\lambda}=(L-\lambda I)^{-1}$, respectively. Also, we denote the eigenvalues of the operators $L_{0}$ and $L$ by  $\>\{\mu_{n}\}_{1}^{\infty} \>$  and $\>\{\lambda_{n}\}_{1}^{\infty}\>$ satisfying the inequalities
$\>\mu_{1}\leq\mu_{2}\leq\cdots\leq\mu_{n}\leq\cdots \>$ and $\>\lambda_{1}\leq\lambda_{2}\leq\cdots\leq\lambda_{n}\leq\cdots\>$.

\noindent If $\>\gamma_{j}\sim aj^{\alpha}\quad(a>0,0<\alpha<\infty)\>$ as $\>j\rightarrow\infty,\>$ then
$$\mu_{n},\lambda_{n}\sim d_{1}n^{\frac{2r\alpha}{2r+\alpha}}, \eqno(1.3)$$

\noindent as $n\rightarrow \infty$ \cite{Yo}. Here,  $\>d_{1}\>$ is a positive constant. By using the asymptotic formula $(1.3)$, there exists a subsequence ${n_{p}}$ of positive integers such that
$$\mu_{q}-\mu_{n_{p}}\geq d_{2}\Bigl(q^{\frac{2r\alpha}{2r+\alpha}}-n_{p}^{\frac{2r\alpha}{2r+\alpha}}\Bigr),\quad
(q=n_{p}+1,n_{p}+2,\ldots) \qquad \eqno(1.4)$$

\noindent where $\>d_{2}\>$ is a positive constant.

\noindent In the work \cite{Pr}, the formula in the form
$$\lim\limits_{p\rightarrow \infty} \sum\limits_{q=1}^{n_{p}}\Biggl[\lambda_{q}-\mu_{q}-\frac{1}{\pi}\int\limits_{0}^{\pi}
\Bigl(Q(x)\varphi_{j_{q}},\varphi_{j_{q}}\Bigr)dx\Biggr]=\frac{1}{4}\Bigl(trQ(0)-trQ(\pi)\Bigr)$$

\noindent is obtained for the  first regularized trace of the operator $L$. In this present work, we find a formula in the form
$$\lim\limits_{p\rightarrow \infty} \sum\limits_{q=1}^{n_{p}}\Biggl(\lambda_{q}^{2}- \mu_{q}^{2}-2\sum\limits_{s=2}^{m}(-1)^{s}s^{-1}
Res_{\lambda=\mu_{q}}tr\bigl[\lambda(Q R_{\lambda}^{0})^{s}\bigr]-\frac{2\mu_{q}}{\pi}\int\limits_{0}^{\pi}
\Bigl(Q(x)\varphi_{j_{q}},\varphi_{j_{q}}\Bigr)dx\Biggr)$$
$$=(-1)^{r}2^{-1-2r}\Bigl[trQ^{(2r)}(0)-trQ^{(2r)}(\pi)\Bigr]+\frac{1}{2}\Bigl[trAQ(0)-trAQ(\pi)\Bigr].\eqno (1.5)$$

\noindent The left hand side of equality (1.5) is called the second regularized trace of the differential operator $L$.

\section{\bf{Main Results}}
\noindent The main purpose of this section is to obtain the second trace formula for the operator $L$. Now, we find the relations between
resolvents and eigenvalues of the operators $L_{0}$ and $L$.

\noindent If $\alpha>\frac{2r}{2r-1}$ and $\lambda\neq\lambda_{q},\mu_{q}$ $(q=1,2,...)$, then by (1.3), $R_{\lambda}^{0}$ and $R_{\lambda}$
are trace class operators. Hence
$$tr(R_{\lambda}-R_{\lambda}^{0})=trR_{\lambda}-trR_{\lambda}^{0}=
\sum\limits_{q=1}^{\infty}\Bigl(\frac{1}{\lambda_{q}-\lambda}-\frac{1}{\mu_{q}-\mu}\Bigr).$$

\noindent  If this equality is multiply with $\>\frac{\lambda^{2}}{2\pi i}\>$ and integrated on the circle

\noindent $\>|\lambda|=b_{p}=\frac{1}{2}(\mu_{n_{p}}+\mu_{n_{p}+1})\>$, then we have following equality
$$\frac{1}{2\pi i}\int\limits_{|\lambda|=b_{p}}\lambda^{2}tr\Bigl(R_{\lambda}-R_{\lambda}^{0}\Bigr)d\lambda=\frac{1}{2\pi i} \int\limits_{|\lambda|=b_{p}} \sum\limits_{q=1}^{\infty}(\frac{\lambda^{2}}{\lambda_{q}-\lambda})d\lambda -\frac{1}{2\pi
i}\int\limits_{|\lambda|=b_{p}} \sum\limits_{q=1}^{\infty}(\frac{\lambda^{2}}{\mu_{q}-\lambda})d\lambda.\eqno(2.1)$$

\noindent We can see that for the large values of $p$,
$$\{\lambda_{q},\mu_{q}\}_{1}^{n_{p}}\subset B(0,b_{p})=\{\lambda:|\lambda|< b_{p}\}$$
$$\lambda_{q},\mu_{q} \notin B[0,b_{p}] =\{\lambda:|\lambda|\leq b_{p}\}\quad (q\geq n_{p}+1).$$

\noindent Therefore by (2.1), we have
$$\sum\limits_{q=1}^{n_{p}}\Bigl(\lambda^{2}_{q}-\mu^{2}_{q}\Bigr)=-{\frac{1}{2\pi i}}\int\limits_{|\lambda|=b_{p}}\lambda^{2}
tr\Bigl(R_{\lambda}-R_{\lambda}^{0}\Bigr) d\lambda.\eqno(2.2)$$

\noindent This is well known formula for the resolvents of the operators $L_{0}$ and $L$:
$$R_{\lambda}=R_{\lambda}^{0}-R_{\lambda}QR_{\lambda}^{0}\qquad (\lambda\in \rho(L)\cap\rho(L_{0})).$$

\noindent By using the last formula, we obtain
$$R_{\lambda}-R_{\lambda}^{0}=\sum\limits_{s=1}^{m}(-1)^{s}R_{\lambda}^{0}(QR_{\lambda}^{0})^{s}+(-1)^{m+1}R_{\lambda}(QR_{\lambda}^{0})^{m+1},$$
\noindent for every positive integer $m$. By (2.2) and the last equality, we have

$$\sum\limits_{q=1}^{n_{p}}\Bigl(\lambda^{2}_{q}-\mu^{2}_{q}\Bigr)={\frac{1}{2\pi i}}\int\limits_{|\lambda|=b_{p}}\lambda^{2}
tr\Bigl(\sum\limits_{s=1}^{m}(-1)^{s+1}R_{\lambda}^{0}(QR_{\lambda}^{0})^{s}+(-1)^{m}R_{\lambda}(QR_{\lambda}^{0})^{m+1}\Bigr) d\lambda$$
\noindent or
$$\sum\limits_{q=1}^{n_{p}}\Bigl(\lambda^{2}_{q}-\mu^{2}_{q}\Bigr)=\sum\limits_{s=1}^{m}D_{ps}+D_{p}^{(m)}.\qquad \qquad \qquad  \eqno(2.3)$$
\noindent Here,
$$D_{ps}={\frac{(-1)^{s+1}}{2\pi i}}\int\limits_{|\lambda|=b_{p}}\lambda^{2} tr\Bigl(R_{\lambda}^{0}(QR_{\lambda}^{0})^{s}\Bigr)d\lambda,
\quad(s=1,2,...) \qquad\eqno(2.4)$$
$$D^{(m)}_{p}={\frac{(-1)^{m}}{2\pi i}}\int\limits_{|\lambda|=b_{p}}\lambda^{2} tr\Bigl(R_{\lambda}(QR_{\lambda}^{0})^{m+1}\Bigr) d\lambda.
\qquad\eqno(2.5)$$

\begin{theorem} If  $\gamma_{j}\sim aj^{\alpha}\quad (0<a, \alpha>\frac{2r}{2r-1})\>$  as $\>j\rightarrow\infty \>$,  then
$$D_{ps}={\frac{(-1)^{s}}{\pi is}}\int\limits_{|\lambda|=b_{p}}\lambda tr\Bigl((QR_{\lambda}^{0})^{s}\Bigr)d\lambda \quad(s=1,2,...).$$

\end{theorem}

\begin{theorem} If the operator function $Q(x)$ satisfies the conditions (Q1) and (Q2), then the series
$$\sum\limits_{k=0}^{\infty}\sum\limits_{j=1}^{\infty}[(k+\frac{1}{2})^{2r}+\gamma_{j}]\int\limits_{0}^{\pi}(Q(x)\varphi_{j,}\varphi_{j})\cos((2k+1)x)
dx$$

\noindent is absolute convergence.

\end{theorem}

\noindent We are at the position to give the main result:

\begin{theorem} If the operator function $\>Q(x)\>$ satisfies the conditions $(Q1)\>$, $\>(Q2)\>$, and  $\gamma_{j}\sim aj^{\alpha} \>$ as
$\>j\rightarrow\infty\quad (a>0,\quad \frac{2r}{2r-1}<\alpha)\>$, then we have

$$\lim_{p\rightarrow\infty}\sum\limits_{q=1}^{n_{p}}\Biggl(\lambda_{q}^{2}-\mu_{q}^{2}-2\sum\limits_{s=2}^{m}(-1)^{s}s^{-1}Res_{\lambda=\mu_{q}}tr\Bigl(\lambda
(QR_{\lambda}^{0})^{s}\Bigr)\Biggr)-\frac{2\mu_{q}}{\pi}\int\limits_{0}^{\pi}(Q(x)\varphi_{j_{q}},\varphi_{j_{q}})dx$$
$$=(-1)^{r}2^{-1-2r}\Bigl[trQ^{(2r)}(0)-trQ^{(2r)}(\pi)\Bigr]+\frac{1}{2}\Bigl[ tr A Q(0)- tr A Q(\pi)\Bigr],\qquad\qquad$$

\noindent where $\>m=\Bigl[\big|\frac{2r\alpha+6r+3\alpha}{2r\alpha -2r-\alpha}\big|\Bigr]$.
\end{theorem}

\section{\bf {Proofs} }

\noindent {\bf Proof of Theorem 1.} We can show that the operator function $(QR_{\lambda}^{0})^{s}$  is analytic with respect to the norm in
the space $\sigma_{1}(H_{1})$ in the region $\rho (L_{0})$ and
$$ tr\Bigl([(QR_{\lambda}^{0})^{s}]'\Bigr)=str\Bigl((QR_{\lambda}^{0})'(QR_{\lambda}^{0})^{s-1}\Bigr),\quad(QR_{\lambda}^{0})'=
Q(R_{\lambda}^{0})^{2}.$$

\noindent Therefore, we have
$$ tr\Bigl([(QR_{\lambda}^{0})^{s}]'\Bigr)=str\Bigl(R_{\lambda}^{0}(QR_{\lambda}^{0})^{s}\Bigr).$$

\noindent From (2.4) and the last equality we obtain
$$ D_{ps}= \frac{(-1)^{s+1}}{2\pi i s}=\int\limits_{|\lambda|=b_{p}}\lambda ^{2}tr\Bigl(\Bigl[(QR_{\lambda}^{0})^{s}\Bigr]'\Bigr)d\lambda.$$

\noindent We can also write the last formula in the following form:
\begin{eqnarray*}
D_{ps}&= &\frac{(-1)^{s+1}}{2\pi i s}\int\limits_{|\lambda|=b_{p}}\Bigl[tr\Bigl(\Bigl[\lambda
^{2}(QR_{\lambda}^{0})^{s}\Bigr]'\Bigr)-2\lambda(QR_{\lambda}^{0})^{s}\Bigr]d\lambda\\
\\
&= &\frac{(-1)^{s}}{\pi i s}\int\limits_{|\lambda|=b_{p}}\lambda tr\Bigl((QR_{\lambda}^{0})^{s}\Bigr)d\lambda+\frac{(-1)^{s+1}}{2\pi i
s}\int\limits_{|\lambda|=b_{p}} tr\Bigl([\lambda^{2}(QR_{\lambda}^{0})^{s}]'\Bigr)d\lambda.\qquad\quad \qquad\qquad \qquad (3.1)
\end{eqnarray*}

\noindent We can see
$$\int\limits_{|\lambda|=b_{p}} tr\Bigl(\Bigl[\lambda^{2}(QR_{\lambda}^{0})^{s}\Bigr]'
\Bigr)d\lambda=\int\limits_{|\lambda|=b_{p}}\Bigl[tr\Bigl(\lambda^{2}(QR_{\lambda}^{0})^{s}\Bigr)\Bigr]'d\lambda.\quad \qquad\eqno(3.2)$$

\noindent  We write the right hand side of above equality in the following way:
$$\int\limits_{|\lambda|=b_{p}}\Bigl[tr\Bigl(\lambda^{2}(QR_{\lambda}^{0})^{s}\Bigr)\Bigr]'d\lambda=\int\limits_{\scriptstyle |\lambda|=b_{p}\atop \scriptstyle Im\lambda\geq 0}
 \Bigl[tr\Bigl(\lambda^{2}(QR_{\lambda}^{0})^{s}\Bigr)\Bigr]'d\lambda+\int\limits_{\scriptstyle |\lambda|=b_{p}\atop \scriptstyle
 Im\lambda\leq 0} \Bigl[tr\Bigl(\lambda^{2}(QR_{\lambda}^{0})^{s}\Bigr)\Bigr]'d\lambda. \quad \qquad \eqno(3.3)$$

\noindent Let $ \varepsilon_{0}$ be a positive number satisfying the inequality $ b_{p}+\varepsilon_{0}<\mu_{n_{p}+1}.$

\noindent We consider that the function $tr\Bigl(\lambda^{2}(QR_{\lambda}^{0})^{s}\Bigr)$ is analytic in the simply connected regions:
$$ G_{1}=\{\lambda:b_{p}- \varepsilon_{0}<|\lambda|<b_{p} +\varepsilon_{0},\quad Im \lambda>-\varepsilon_{0}\},$$
$$ G_{2}=\{\lambda:b_{p}- \varepsilon_{0}<|\lambda|<b_{p} +\varepsilon_{0},\quad Im \lambda<\varepsilon_{0}\}$$

\noindent and
 $$\{\lambda:|\lambda|=b_{p},\quad Im \lambda\geq 0\}\subset G_{1},$$
 $$\{\lambda:|\lambda|=b_{p},\quad Im \lambda\leq 0\}\subset G_{2}.$$

\noindent By using the Leibnitz Formula and (3.3), we get
$$\int\limits_{|\lambda|=b_{p}}\Big\{tr\Bigl(\lambda^{2}(QR_{\lambda}^{0})^{s}\Bigr)\Big\}'d\lambda\qquad\qquad\qquad\qquad\qquad\qquad\qquad\qquad$$
$$= tr\Bigl(b_{p}^{2}(QR_{-b_{p}}^{0})^{s}\Bigr) -
tr\Bigl(b_{p}^{2}(QR_{b_{p}}^{0})^{s}\Bigr)+tr\Bigl(b_{p}^{2}(QR_{b_{p}}^{0})^{s}\Bigr)-tr\Bigl(b_{p}^{2}(QR_{-b_{p}}^{0})^{s}\Bigr)=0.\eqno(3.4)$$
\noindent From (3.1), (3.2) and (3.4), we have
$$D_{ps}={\frac{(-1)^{s}}{\pi is}}\int\limits_{|\lambda|=b_{p}}\lambda tr\Bigl((QR_{\lambda}^{0})^{s}\Bigr)d\lambda .\qquad\qquad\qquad\qquad\qquad\qquad\qquad\qquad\Box$$

\noindent {\bf Proof of Theorem 2.} Let $h_{j}(x)=(Q(x)\varphi_{j,}\varphi_{j})$. Using the integration by parts formula and the condition (Q1), we get

\begin{eqnarray*}
& &\int\limits_{0}^{\pi}h_{j}(x)\cos((2k+1)x)dx=\int\limits_{0}^{\pi}h_{j}(x)\Bigl(\frac{1}{2k+1}\sin(2k+1)x\Bigr)'dx\\
\\
&=&\frac{1}{2k+1}\Bigl[h_{j}(x)\sin((2k+1)x)\Bigl|^{\pi}_{0}-\int\limits_{0}^{\pi}h'_{j}(x)\sin((2k+1)x)dx\Bigr]\\
\\
&=&\frac{1}{(2k+1)^2}\int\limits_{0}^{\pi}h'_{j}(x)\bigl(\cos(2k+1)x\bigr)'dx\\
\\
&=&\frac{1}{(2k+1)^2}\Bigl[h'_{j}(x)\cos((2k+1)x)\Bigl|^{\pi}_{0}-\int\limits_{0}^{\pi}h''_{j}(x)\cos((2k+1)x)dx\Bigr]\\
\\
&=&-\frac{1}{(2k+1)^3}\int\limits_{0}^{\pi}h''_{j}(x)\bigl(\sin(2k+1)x\bigr)'dx\\
\\
&=&-\frac{1}{(2k+1)^3}\Bigl[h_{j}''(x)\sin((2k+1)x)\Bigl|^{\pi}_{0}-\int\limits_{0}^{\pi}h''_{j}(x)\sin((2k+1)x)dx\Bigr]\\
\\
&=&....=\frac{(-1)^{r+1}}{(2k+1)^{2r+2}}\int\limits_{0}^{\pi}h^{(2r+2)}_{j}(x)\cos((2k+1)x)dx.\quad\qquad\qquad\qquad\qquad\qquad\qquad\qquad\qquad(3.5)
\end{eqnarray*}

\noindent By (3.5), we find
\begin{eqnarray*}
&  &\sum\limits_{k=0}^{\infty}\sum\limits_{j=1}^{\infty}\Bigl|[(k+\frac{1}{2})^{2r}+\gamma_{j}\Bigr]\int\limits_{0}^{\pi}h_{j}(x)\cos((2k+1)x)
dx\Bigr|\\
\\
&\leq &\sum\limits_{k=0}^{\infty}\sum\limits_{j=1}^{\infty}(2k+1)^{-2}\int\limits_{0}^{\pi}\Biggl[\Big|h^{(2r+2)}_{j}(x)\Big|+\gamma_{j}|h''_{j}(x)|\Biggr]dx\\
\\
&= & \sum\limits_{j=1}^{\infty}\int\limits_{0}^{\pi}\Biggl[ \Big|(Q^{(2r+2)}(x)\varphi_{j,}\varphi_{j})\Big|+|(A
Q''(x)\varphi_{j,}\varphi_{j})|\Biggr]dx \sum\limits_{k=0}^{\infty}(2k+1)^{-2}\\
\\
&\leq & Const.\int\limits_{0}^{\pi}\Biggl[\sum\limits_{j=1}^{\infty}
\Big|(Q^{(2r+2)}(x)\varphi_{j,}\varphi_{j})\Big|+\sum\limits_{j=1}^{\infty}|(A
Q''(x)\varphi_{j,}\varphi_{j})|\Biggr]dx
\end{eqnarray*}
$$\leq Const.\int\limits_{0}^{\pi}\Biggl[\Big\|Q^{(2r+2)}(x)\Big\|_{\sigma_{1}(H)}+\|A
Q''(x)\|_{\sigma_{1}(H)}\Biggr]dx.\qquad\qquad\eqno(3.6)$$

\noindent Since the functions $\| Q^{(2r+2)}(x)\|_{\sigma_{1}(H)}$ and $\|A Q''(x)\|_{\sigma_{1}(H)}$ in (3.6) are measurable and bounded  in
the interval $[0,\pi]$, we get
$$\sum\limits_{k=0}^{\infty}\sum\limits_{j=1}^{\infty}\Bigg|[(k+\frac{1}{2})^{2r}+\gamma_{j}]\int\limits_{0}^{\pi}(Q(x)\varphi_{j,}\varphi_{j})\cos((2k+1)x)
dx\Bigg|<\infty.\qquad\qquad\qquad\quad\quad\Box$$

\noindent Let $\{\psi_{q}\}_{1}^{\infty}$ be the orthornormal eigenvectors system  corresponding to eigenvalues $\{\mu_{q}\}_{1}^{\infty}$ of
the operator $L_{0}$, respectively. Since the orthornormal eigenvectors  corresponding to eigenvalues

\noindent $(k+\frac{1}{2})^{2r}+\gamma_{j}\qquad
(k=0,1,2,\cdots;j=1,2,\cdots)$ of the operator $L_{0}$ are $\sqrt{\frac{2}{\pi}}\cos((k+\frac{1}{2})x)\varphi_{j}$, respectively,
$$ \mu_{q}= (k+\frac{1}{2})^{2r}+\gamma_{j_{q}}\qquad (q=1,2,\cdots)$$

\noindent and
$$\psi_{q}(x)= \sqrt{\frac{2}{\pi}}\cos((k_{q}+\frac{1}{2})x)\varphi_{j_{q}}. \qquad\quad\qquad\eqno(3.7)$$

\noindent We prove the main theorem of the paper.

\noindent \textbf{Proof of Theorem 3.}

\noindent By using the Theorem 1, one can write $\>D_{ps}\>$ as follows:
\begin{eqnarray*}
D_{ps}&= &2(-1)^{s}s^{-1}\frac{1}{2\pi i} \int\limits_{|\lambda|= b_{p}}tr\Bigl(\lambda(QR_{\lambda}^{0})^{s}\Bigr)d\lambda\\
\\
&= &2(-1)^{s}s^{-1}\sum\limits_{q=1}^{n_{p}}Res_{\lambda=\mu_{q}}tr\Bigl(\lambda(QR_{\lambda}^{0})^{s}\Bigr).
\end{eqnarray*}

\noindent By using last formula, we can rewrite (2.3) as follows:
$$\sum\limits_{q=1}^{n_{p}}\Biggl(\lambda_{q}^{2}-\mu_{q}^{2}-2\sum\limits_{s=2}^{m}(-1)^{s}s^{-1}Res_{\lambda=\mu_{q}}tr\Bigl(\lambda
(QR_{\lambda}^{0})^{s}\Bigr)\Biggr)=D_{p1}+D_{p}^{(m)},\qquad\eqno(3.8)$$
$$D_{p1}=\frac{-1}{\pi i} \int\limits_{|\lambda|= b_{p}}\lambda tr(QR_{\lambda}^{0})d\lambda.\qquad\eqno(3.9)$$

\noindent Since $(QR_{\lambda}^{0})$ is a nuclear operator for every $\lambda \in \rho(L_{0})$ and  $\{\psi_{q}\}_{1}^{\infty}$ is an
orthonormal basis in the space $H_{1}$, we have
$$tr(QR_{\lambda}^{0})=\sum\limits_{q=1}^{\infty}(QR_{\lambda}^{0}\psi_{q},\psi_{q})_{{H}_{1}}.$$
\noindent Here, $ R_{\lambda}^{0} \psi_{q}=(\mu_{q}-\lambda I)^{-1}\psi_{q}.$
\noindent If we substitute the last two equalities into (3.9), then we get
\begin{eqnarray*}
D_{p1}&= &\frac{-1}{\pi i} \int\limits_{|\lambda|= b_{p}}\lambda
\sum\limits_{q=1}^{\infty}(QR_{\lambda}^{0}\psi_{q},\psi_{q})_{{H}_{1}}d\lambda\\
\\
&= & \frac{-1}{\pi i} \int\limits_{|\lambda|= b_{p}}\lambda
\sum\limits_{q=1}^{\infty}\frac{1}{\mu_{q}-\lambda}(Q\psi_{q},\psi_{q})_{{H}_{1}}d\lambda\\
\\
&= & \frac{1}{\pi i}\sum\limits_{q=1}^{\infty}(Q\psi_{q},\psi_{q})_{{H}_{1}}\int\limits_{|\lambda|=
b_{p}}\frac{\lambda}{\lambda-\mu_{q}}d\lambda
\end{eqnarray*}

\noindent By using the  Cauchy Integral Formula
\[ \frac{1}{2\pi i}\int\limits_{|\lambda|= b_{p}}\frac{\lambda}{\lambda-\mu_{q}}d\lambda=   \left\{\begin{array}{cl} \mu_{q} &, if \qquad
q\leq n_{p}\\0 &, if \qquad q>n_{p} \end{array} \right. \]

\noindent and by (3.7), we obtain
\begin{eqnarray*}
D_{p1}&= &2 \sum\limits_{q=1}^{n_{p}}\mu_{q}(Q\psi_{q},\psi_{q})_{{H}_{1}}\\
\\
&= &2 \sum\limits_{q=1}^{n_{p}}\mu_{q}\int\limits_{0}^{\pi}(Q(x)\psi_{q}(x),\psi_{q}(x))dx\\
\\
&= & 2 \sum\limits_{q=1}^{n_{p}}\mu_{q}\int\limits_{0}^{\pi}\Bigl(Q(x)\sqrt{\frac{2}{\pi}}\cos
(k_{q}+\frac{1}{2})x\varphi_{j_{q}},\sqrt{\frac{2}{\pi}}\cos (k_{q}+\frac{1}{2})x\varphi_{j_{q}}\Bigr)dx\\
\\
&= &  2
\sum\limits_{q=1}^{n_{p}}\mu_{q}\frac{2}{\pi}\int\limits_{0}^{\pi}\cos^{2}(k_{q}+\frac{1}{2})x(Q(x)\varphi_{j_{q}},\varphi_{j_{q}})dx\\
\\
&= & \sum\limits_{q=1}^{n_{p}}\mu_{q}\frac{2}{\pi}\int\limits_{0}^{\pi}(1+\cos(2k_{q}+1)x)(Q(x)\varphi_{j_{q}},\varphi_{j_{q}})dx\\
\\
&= &\frac{2}{\pi} \sum\limits_{q=1}^{n_{p}}\mu_{q}\int\limits_{0}^{\pi}\cos(2k_{q}+1)x(Q(x)\varphi_{j_{q}},\varphi_{j_{q}})dx+\frac{2}{\pi}
\sum\limits_{q=1}^{n_{p}}\mu_{q}\int\limits_{0}^{\pi}(Q(x)\varphi_{j_{q}},\varphi_{j_{q}})dx.\quad(3.10)
\end{eqnarray*}

\noindent We substitude (3.10) in  (3.8):
\begin{eqnarray*}
& &\sum\limits_{q=1}^{n_{p}}\Biggl(\lambda_{q}^{2}-\mu_{q}^{2}-2\sum\limits_{s=2}^{m}(-1)^{s}s^{-1}Res_{\lambda=\mu_{q}}tr\Bigl[\lambda
(QR_{\lambda}^{0})^{s}\Bigr]-\frac{2\mu_{q}}{\pi}\int\limits_{0}^{\pi}h_{j_{q}}(x)dx\Biggr)\\
\\
&\qquad= &\frac{2}{\pi} \sum\limits_{q=1}^{n_{p}}\mu_{q}\int\limits_{0}^{\pi}h_{j_{q}}(x) \cos(2k_{q}+1)x
dx+D_{p}^{(m)}.\qquad\qquad\qquad\qquad\qquad\qquad\qquad(3.11)
\end{eqnarray*}

\noindent If we use Theorem 2, then we know that
\begin{eqnarray*}
& &\frac{2}{\pi}\lim_{p\rightarrow\infty}\sum\limits_{q=1}^{n_{p}}\mu_{q}\int\limits_{0}^{\pi}h_{j_{q}}(x) \cos(2k_{q}+1)x dx\\
\\
&\qquad= &\frac{2}{\pi}
\sum\limits_{k=0}^{\infty}\sum\limits_{j=1}^{\infty}\Bigl((k+\frac{1}{2})^{2r}+\gamma_{j}\Bigr)\int\limits_{0}^{\pi}h_{j}(x) \cos(2k+1)x dx.
\qquad\qquad\qquad\qquad\qquad(3.12)
\end{eqnarray*}

\noindent If we substitude (3.5) in (3.12), then we get

\begin{eqnarray*}
& &\frac{2}{\pi} \sum\limits_{k=0}^{\infty}\sum\limits_{j=1}^{\infty}\Bigl[(k+\frac{1}{2})^{2r}+\gamma_{j}\Bigr]\int\limits_{0}^{\pi}h_{j}(x)
\cos(2k+1)x dx\\
\\
&=
&\sum\limits_{k=0}^{\infty}\sum\limits_{j=1}^{\infty}\frac{2}{\pi}\int\limits_{0}^{\pi}\Bigl[(-\frac{1}{4})^{r}h_{j}^{(2r)}(x)+\gamma_{j}(x)h_{j}(x)\Bigr]
\cos(2k+1)x dx\\
\\
&=
&\frac{1}{\pi}\sum\limits_{j=1}^{\infty}\sum\limits_{k=0}^{\infty}\Biggl(\int\limits_{0}^{\pi}\Bigl[(-\frac{1}{4})^{r}h_{j}^{(2r)}(x)+\gamma_{j}(x)h_{j}(x)\Bigr]
\cos(k x) dx\\
\\
&- &(-1)^{k}\int\limits_{0}^{\pi}\Bigl[(-\frac{1}{4})^{r}h_{j}^{(2r)}(x)+\gamma_{j}(x)h_{j}(x)\Bigr] \cos(k x)dx \Biggr)\\
\\
&=
&\frac{1}{2}\sum\limits_{j=1}^{\infty}\Biggl\{\sum\limits_{k=0}^{\infty}M_{k}\int\limits_{0}^{\pi}\Bigl[(-\frac{1}{4})^{r}h_{j}^{(2r)}(x)+\gamma_{j}(x)h_{j}(x)\Bigr]
\cos(k x) dx.\cos(k0)\\
\\
&- &\sum\limits_{k=0}^{\infty}M_{k}\int\limits_{0}^{\pi}\Bigl[(-\frac{1}{4})^{r}h_{j}^{(2r)}(x)+\gamma_{j}(x)h_{j}(x)\Bigr] \cos(k x)dx
\cos(k\pi) \Biggr\},
\end{eqnarray*}

\noindent here, \[ M_{k}=  \left\{ \begin{array}{cl} \pi^{-1} &, if \qquad k=0\\2\pi^{-1} &, if \qquad k=1,2,\ldots \end{array} \right. \]

\noindent The sums according to the $k$ on the right hand side of the last relation are the values at $0$ and $\pi$ of the Fourier Series of
the function $-(\frac{1}{4})^{r}h_{j}^{(2r)}(x)+\gamma_{j}h_{j}(x) $  according to the  functions $\{\cos kx\}^{\infty}_{k=0} $ on the
interval $[0,\pi]$.

\noindent Therefore,
\begin{eqnarray*}
& &\frac{2}{\pi} \sum\limits_{k=0}^{\infty}\sum\limits_{j=1}^{\infty}\Bigl[(k+\frac{1}{2})^{2r}+\gamma_{j}\Bigr]\int\limits_{0}^{\pi}h_{j}(x)
\cos(2k+1)x dx\\
\\
&= &\frac{1}{2}\sum\limits_{j=1}^{\infty}\Bigl[(-\frac{1}{4})^{r}(h_{j}^{(2r)}(0)-h_{j}^{(2r)}(\pi))+\gamma_{j}(h_{j}(0)-h_{j}(\pi))\Bigr] \\
\\
&= &(-1)^{r}2^{-1-2r}\Bigl[trQ^{(2r)}(0)-trQ^{(2r)}(\pi)\Bigr]+\frac{1}{2}\Bigl[ tr A Q(0)- tr A Q(\pi)\Bigr]\qquad\qquad\qquad \qquad\qquad(3.13)
\end{eqnarray*}

\noindent From (3.12) and (3.13), we obtain
\begin{eqnarray*}
& &\frac{2}{\pi}\lim_{p\rightarrow\infty}\sum\limits_{q=1}^{n_{p}}\mu_{q}\int\limits_{0}^{\pi}h_{j_{q}}(x) \cos(2k_{q}+1)x dx\\
\\
&= &(-1)^{r}2^{-1-2r}\Bigl[trQ^{(2r)}(0)-trQ^{(2r)}(\pi)\Bigr]+\frac{1}{2}\Bigl[ tr A Q(0)- tr A Q(\pi)\Bigr]\qquad\qquad\qquad\qquad\qquad(3.14)
\end{eqnarray*}

\noindent Let us estimate of $D_{p}^{(m)}$ for the large value of p. By using (2.5), we get
\begin{eqnarray*}
\Big|D_{p}^{(m)}\Big|&\leq &\int\limits_{|\lambda|= b_{p}}|\lambda|^{2}\Big|tr\bigl(R_{\lambda}(QR_{\lambda}^{0})^{m+1}\bigr)\Big|
|d\lambda|\\
\\
&\leq & b_{p}^{2}\int\limits_{|\lambda|= b_{p}}\Big\|R_{\lambda}(QR_{\lambda}^{0})^{m+1}\Big\|_{\sigma_{1}(H_{1})} |d\lambda|\\
\\
&\leq & b_{p}^{2}\int\limits_{|\lambda|= b_{p}}\Big\|R_{\lambda}\Big\|_{1}\Big\|(QR_{\lambda}^{0})^{m+1}\Big\|_{\sigma_{1}(H_{1})}
|d\lambda|\\
\\
&\leq & b_{p}^{2}\int\limits_{|\lambda|= b_{p}}
\Big\|R_{\lambda}\Big\|_{1}\Big\|(QR_{\lambda}^{0})^{m}\Big\|_{1}\Big\|QR_{\lambda}^{0}\Big\|_{\sigma_{1}(H_{1})} |d\lambda|\\
\\
&\leq & b_{p}^{2}\int\limits_{|\lambda|= b_{p}} \Big\|R_{\lambda}\Big\|_{1}\Big\|Q\Big\|_{1}\Big\|R_{\lambda}^{0}\Big\|^{m}_{1}\Big\|Q\Big\|
_{1}\Big\| R_{\lambda}^{0}\Big\|_{\sigma_{1}(H_{1})} |d\lambda|.\quad\qquad\qquad\qquad\qquad\qquad\quad\quad\quad(3.15)
\end{eqnarray*}

\noindent One can prove the following inequalities similarly in work \cite{Yo}:

$$\Big\| R_{\lambda}^{0}\Big\|_{\sigma_{1}(H_{1})}\leq const. {n_{p}}^{1-\delta},$$

$$\Big\| R_{\lambda}\Big\|_{1}\leq const. {n_{p}}^{-\delta} \quad(\delta=\frac{2r\alpha}{2r+\alpha}-1).$$

\noindent From last two inequalities and (3.15), we obtain
$$\Big|D_{p}^{(m)}\Big|\leq const. b_{p}^{3}{n_{p}}^{-\delta m-2\delta+1}.\qquad \eqno(3.16)$$

\noindent For large values of p
$$\>b_{p}=2^{-1}\bigl(\mu_{n_{p}}+\mu_{n_{p}+1}\bigr)\leq const.n_{p}^{1+\delta}.\eqno(3.17)$$

\noindent From (3.16) and (3.17), we obtain
$$\Big|D_{p}^{(m)}\Big|\leq const. \mu_{p}^{4-(m-1)\delta}.\qquad $$

\noindent Therefore, for $\>m=\Bigl[\big|\frac{2r\alpha+6r+3\alpha}{2r\alpha -2r-\alpha}\big|\Bigr]+1$, we find
$$\lim_{p\rightarrow\infty}D_{p}^{(m)}=0 .\qquad\qquad \eqno(3.18)$$

\noindent From $\>(3.11), (3.14)\>$ and $\>(3.18)$, we find the following formula for second regularized trace formula of the operator
$\>L$
$$\lim_{p\rightarrow\infty}\sum\limits_{q=1}^{n_{p}}\Biggl(\lambda_{q}^{2}-\mu_{q}^{2}-2\sum\limits_{s=2}^{m}(-1)^{s}s^{-1}Res_{\lambda=\mu_{q}}tr\Bigl(\lambda
(QR_{\lambda}^{0})^{s}\Bigr)\Biggr)-\frac{2\mu_{q}}{\pi}\int\limits_{0}^{\pi}(Q(x)\varphi_{j_{q}},\varphi_{j_{q}})dx$$
$$=(-1)^{r}2^{-1-2r}\Bigl[trQ^{(2r)}(0)-trQ^{(2r)}(\pi)\Bigr]+\frac{1}{2}\Bigl[ tr A Q(0)- tr A Q(\pi)\Bigr].\qquad\qquad\Box$$

\end{document}